\documentclass[10pt,a4paper]{article}
\usepackage{indentfirst}
\usepackage{amsmath,amssymb,amsfonts,amsthm,graphics}
\usepackage{epsfig}
\usepackage{makeidx}
\usepackage{mathrsfs}
\usepackage{color}
\usepackage{cite}
\usepackage{dsfont}
\usepackage{mathrsfs}
\usepackage{algorithmic}
\usepackage{algorithm}
\usepackage{graphicx}
\usepackage{geometry}
\usepackage{epstopdf}
\usepackage{setspace}
\usepackage{caption}
\usepackage{subfig}
\usepackage[section]{placeins}
\usepackage{float}
\newtheorem{theorem}{Theorem}[section]

\newtheorem{lemma}[theorem]{Lemma}

\newtheorem{remark}[theorem]{Remark}

\newtheorem{conjecture}{Conjecture}[section]
\newtheorem{problem}[conjecture]{Problem}

\begin{document}
\title{\Large\bf Diameter two orientability of mixed graphs}
\author{Hengzhe Li$^{a}$, Zhiwei Ding$^{a}$, Jianbing Liu$^{b}$\footnote{Corresponding author. Email: jianliu@hartford.edu}, Hong-Jian Lai$^{c}$ \\
\small $^{a}$College of Mathematics and Information Science,\\
		\small Henan Normal University, Xinxiang 453007, P.R. China\\
		\small $^b$ Department of Mathematics, \\
        \small University of Hartford, West Hartford 06117, USA\\
        \small $^c$Department of Mathematics\\
		\small West Virginia University, Morgantown, WV 26506, United States\\
        \small Email: lihengzhe@htu.edu.cn, dzhiweidzw@163.com, jianliu@hartford.edu, hjlai2015@hotmail.com;\\}
\date{}
\maketitle

\begin{abstract}
In 1967, Katona and Szemer\'{e}di showed that
no undirected graph with $n$ vertices and fewer than $\frac{n}{2}\log_2\frac{n}{2}$ edges admits an orientation of diameter two.
In 1978, Chv\'atal and Thomassen revealed the complexity of determining whether an undirected graph can be oriented to achieve a diameter of two, proving it to be NP-complete. This breakthrough has sparked ongoing interest in identifying sufficient conditions for graphs to be oriented with the smallest possible diameter of two---critical for optimizing communication and network flow in larger structures. In 2019, Czabarka, Dankelmann, and Sz\'ekely significantly advanced this field by establishing that the minimum degree threshold for achieving such an orientation in undirected graphs of order $n$ is $\frac{n}{2} + \Theta(\ln n)$. In this paper, we extend this foundational result by determining the minimum degree threshold necessary for realizing an orientation with diameter two in mixed graphs, which contain both undirected and directed edges. Mixed graphs offer a versatile framework, representing an intermediate stage in the orientation process, making our findings a substantial generalization of previous results.

\end{abstract}

\section{Introduction}
In this paper, all graphs we considered are mixed graphs without parallel edges and loops. We refer to book \cite{Bondy} for undefined notation in the following.

Mixed graphs, which contains both undirected and directed edges, serve as a generalization of undirected and directed graphs and can be viewed as an intermediate stage in the graph orientation process.

Let $G$ be a mixed graph, $V(G)$ and $E(G)$ represent the vertex set and edge set of $G$, respectively. For a vertex $x$ in $G$, $N_G^+(x)=\{y:\overrightarrow{xy}\in E(G)\}$, $N_G^-(x)=\{y:\overrightarrow{yx}\in E(G)\}$,
$N_G^-(x)=\{y:xy\in E(G)\}$, and $N_G(x)=N_G^+(x)\cup N_G^0(x)\cup N_G^0(x)$. The {\it degree} ({\it out-degree}, {\it in-degree}, {\it un-degree}, respectively) of $x$ is
$d_G(x)=|N_G(x)|$ ($d_G^+(x)=|N_G^+(x)|$, $d_G^-(x)=|N_G^-(x)|$,
$d_G^0(x)=|N_G^0(x)|$, respectively).
Let $\delta(G)$ and $\Delta(G)$ represent the {\it minimum degree} and {\it maximum degree} of the vertices of $G$.

To measure the ratio of directed edges in a mixed graph, we introduce the concept of the mixed ratio of a mixed graph as follow.
For two rational numbers $c_1$ and $c_2$ such that $0\leq c_1,c_2\le 1$, a mixed graph $G$ is a {\it $(c_1,c_2)$-mixed graph} if $d_G^+(u)\leq c_1d_G(u)$ and $d_G^-(u)\leq c_2d_G(u)$ for each $u\in V(G)$,
and $(c_1,c_2)$ is called the {\it mixed ratio} of $G$. Thus, for each mixed graph $G$, there exist rational numbers $c_1$ and $c_2$ with $0\leq c_1,c_2\le 1$ such that $G$ is a {\it $(c_1,c_2)$-mixed graph}.

In a mixed $G$, an {\it $s$-path} is a sequence of distinct vertices $y_1,y_2,\cdots y_s$ such that $y_iy_{i+1}\in E(G)$ or $\overrightarrow{y_iy_{i+1}}\in E(G)$ for $1\leq i\leq s-1$, a {\it directed $s$-path} is a sequence of distinct vertices $y_1,y_2,\cdots y_s$
such that $\overrightarrow{y_iy_{i+1}}\in E(G)$ for $1\leq i\leq s-1$.
The mixed graph is connected if there is a path from $u$ to $v$ for any two distinct vertices $u,v\in V(G)$. The number of edges in the shortest path from $u$ to $v$ is called the {\it distance} from $u$ to $v$ and denoted by $d_G(u,v)$. The diameter {\it diam(H)} is the largest distance between two vertices.

In a mixed graph $G$, a {\it bridge} is an undirected edge $e$ if $G-e$ is disconnected. If the graph is connected and has no bridge, we call it a {\it birdgeless graph}. Specify a direction for all undirected edges in $G$ is called an {\it orientation} of $G$. The graph $G$ is {\it strongly orientable} if it has a strongly connected orientation.
The {\it oriented diameter} of a mixed graph is the minimum diameter among all its strong orientations, and denoted $\overrightarrow{diam}(G)$.

The well-known Robbin's Theorem \cite{Robbins} states that a graph $G$ has a strong orientation if $G$ is a birdgeless connected graph.
Based on this result, it is natural to find an orientation of a graph with minimum diameter. In 1967, Katona and Szemer\'{e}di showed that
no undirected graph with $n$ vertices and fewer than $\frac{n}{2}\log_2\frac{n}{2}$ edges admits an orientation of diameter two.
In 1978, Chv\'atal and Thomassen \cite{Chvatal} showed that there exists a function $f$
satisfying for every bridgeless graph with diameter $d$ has an oriented diameter at most $f(d)$. Specifically, they proved that $f(2)=6$ and that $\frac{1}{2}d^{2}+d\leq f(d)\leq 2d^{2}+2d$. Obviously, this shows that $8\leq f(3)\leq 24$. Subsequently, Kwok, Liu and West \cite{Kwok} improved the range and proved $9\leq f(3)\leq 11$. Ultimately, in 2022, Wang and Chen \cite{Wang} conclusively determined that $f(3)=9$.

Upper bounds on $\overrightarrow{diam}(G)$ have been established in terms of various graph parameters, including maximum degree \cite{B.Chen, P.Dankelmann}, the minimum degree \cite{Bau,Surmacs} and domination number \cite{Dankelmann,Fomin,Kurz}, and radius \cite{Babu,Li}. The oriented radius and oriented diameter have also been investigated for special graph classes \cite{Cochran,Kumar,Mondal}.

Studying the sufficient conditions for graphs that can be oriented to achieve a diameter of two is interesting, as this is the smallest possible diameter for orientations of graphs with more than one vertex. In 1978, Chv\'atal and Thomassen \cite{Chvatal} established that determining an undirected graph has an orientation of diameter two is $NP$-complete. Let $\delta(n)$ be the smallest value such that every graph of order $n$ and minimum degree at least $\delta(n)$ admits oriented diameter two. In 2019, Czabarka, Dankelmann and Sz\'ekely \cite{Czabarka} showed that $\delta(n)=\frac{n}{2}+\Theta(\ln n)$ for $n\in\mathbb{N}$. Subsequently, Chen and Chang \cite{Chen} considered diameter three orientability of bipartite graphs, and showed that for balanced bipartite graphs of order $n$, the minimum degree threshold for diameter three orientability is $\frac{n}{4}+\Theta(\ln n)$.

In this paper, we establish a similar degree condition for diameter two orientability of mixed graphs by showing that $\delta(n,c_1,c_2)=\frac{n}{2-c_1-c_2}+\Theta(\ln n)$. This result can be viewed as a generalization of the result by Czabarka, Dankelmann and Sz\'ekely, since $\delta(n,0,0)=\frac{n}{2}+\Theta(\ln n)$ when $c_1=c_2=0$.

\section{Main Result}
Firstly, the following theorem states that $\delta(n,c_1,c_2)$ is at most $\frac{n}{2-c_1-c_2}+\Theta(\ln n)$.

\begin{theorem}\label{Upper} Let $c_1$, $c_2$ be two rational numbers with $0\leq c_1,c_2\le 1$, and let $G$ be a $(c_1,c_2)$-mixed graph with $n$ vertices. If
$$\delta(G)\geq\frac{1}{2-c_1-c_2}n+\frac{2}{2-c_1-c_2}
  \frac{\ln n}{\ln \frac{4}{3}},$$
then $\overrightarrow{diam}(G)=2$.
\end{theorem}

\begin{pf} Assume $\delta(G)\ge \frac{n}{2}+h(n)$, where $h(n)$ is a function of $n$. Then, for each $u\in V(G)$, we have
$$d_G^+(u)+d_G^0(u)= d_G(u)-d_G^-(u)\geq (1-c_2)d_G(u),$$
$$d_G^-(u)+d_G^0(u)= d_G(u)-d_G^+(u)\geq (1-c_1)d_G(u).$$
Orient every undirected edge of $G$ randomly and independently with probability $\frac{1}{2}$.
Define the random variable $X_{uv}$ as follows:

$$X_{uv}=\left\{\begin{array}{ll}
1, & \text{there is no directed $2$-path from $u$ to $v$};\\
0, & \text{Otherwise}.
 \end{array}\right.$$

Let $P(X_{uv})$ be the probability of the event $X_{uv}$, let $X=\sum_{u,v\in V(G), u\neq v} X_{uv}$, and let $\xi$ be the expected value of $X$. By the linearity of expectation,  we have
$$\xi=\sum_{u,v\in V(G),u\neq v }P[X_{uv}].$$

For two distinct vertices $u,v\in V(G)$, to estimate $P[X_{uv}=1]$, we need to consider the following two cases.

\noindent\textit{Case 1.} $N_G^+(u)\cap N_G^-(v)\neq \emptyset$.

In this case, we have $X_{uv}=0$ since there is a directed $2$-path
from $u$ to $v$. Hence, $P[X_{uv}=1]=0$.

\noindent\textit{Case 2.} $N_G^+(u)\cap N_G^-(v)=\emptyset$.

Let $B=N_G(u)\cap N_G(v)$. Since $N_G^+(u)\cap N_G^-(v)=\emptyset$, the set $B$ can be divided into the following three subsets.
\begin{align*}
B_1&=N_G^+(u)\cap N_G^0(v),\\
B_2&=N_G^0(u)\cap N_G^-(v),\\
B_3&=N_G^0(u)\cap N_G^0(v).
\end{align*}

By the inclusion-exclusion formula,
\begin{align*}
&|B|\geq (1-c_2)d_G(u)+(1-c_1)d_G(v)-n\ge
(2-c_1-c_2)h(n)-\frac{c_1+c_2}{2}n,\\
&|B_3|\geq (1-c_1-c_2)d_G(u)+(1-c_1-c_2)d_G(v)-n
\ge 2(1-c_1-c_2)h(n)-(c_1+c_2)n.
\end{align*}

Thus, we obtain
\begin{align*}
 P[X_{uv}=1]&=\left(\frac{1}{2}\right)^{|B_1|} \left(\frac{1}{2}\right)^{|B_2|} \left(\frac{3}{4}\right)^{|B_3|}\\
 &= \left(\frac{1}{2}\right)^{|B|} \left(\frac{3}{2}\right)^{|B_3|}\\
 &\leq\left(\frac{1}{2}\right)^{(2-c_1-c_2)h(n)-\frac{c_1+c_2}{2}n}  \left(\frac{3}{2}\right)^{2(1-c_1-c_2)h(n)-(c_1+c_2)n}\\
 &=\left(\frac{3}{4}\right)^{(2-c_1-c_2)h(n)-\frac{c_1+c_2}{2}n}  \left(\frac{3}{2}\right)^{-(c_1+c_2)h(n)-\frac{c_1+c_2}{2}n}\\
 &\leq\left(\frac{3}{4}\right)^{(2-c_1-c_2)h(n)-\frac{c_1+c_2}{2}n}
\end{align*}

Therefore, we have
$$\xi = \sum_{u,v\in V(G),u\neq v }P[X_{uv}]\leq n(n-1)  \left(\frac{3}{4}\right)^{(2-c_1-c_2)h(n)-\frac{c_1+c_2}{2}n}
<n^{2}   \left(\frac{3}{4}\right)^{(2-c_1-c_2)h(n)-\frac{c_1+c_2}{2}n}$$

If $h(n)$ is large enough, then
$$\xi<n^{2}\left(\frac{3}{4}\right)^{(2-c_1-c_2)h(n)-\frac{c_1+c_2}{2}n}
\leq1$$

Therefore, the graph $G$ can be oriented in such a way that for any two vertices $u, v\in V(G)$, there is a directed $2$-path from $u$ to $v$. Through straightforward calculations, one derives that $h(n)\geq\frac{2}{2-c_1-c_2} \frac{n}{2}+\frac{2}{2-c_1-c_2}
 \frac{\ln n}{\ln \frac{4}{3}}$ implies that $\xi<1$. Consequently, the graph $G$ has the oriented diameter two if $ \delta(G)\geq\frac{n}{2}+h(n)\ge\frac{1}{2-c_1-c_2}n+\frac{2}{2-c_1-c_2}
  \frac{\ln n}{\ln \frac{4}{3}}$. The theorem is proved.
\end{pf}$\hfill\square$

Next, we construct an infinite family of $(c_1,c_2)$-mixed graphs $G_m$ with minimum degree $\frac{n}{2-c_1-c_2}+\Theta(\ln n)$, but $\overrightarrow{diam}(G_m)\geq 3$.

\begin{lemma}\label{Lower1} Let $c_1$, $c_2$ be two rational numbers with $0\leq c_1, c_2\le 1$. If $\frac{1-c_1-c_2}{c(2-c_1-c_2)}\geq 1$, then there exists an infinite family of $(c_1,c_2)$-mixed graphs $G_m$ of order $n$ and
$$\delta(G_m)\geq\frac{1}{2-c_1-c_2}n+\frac{2}{2-c_1-c_2}
  \frac{\ln n}{2\ln \frac{27}{4}},$$
whose oriented diameter is at least three.
\end{lemma}

\begin{pf}
Let $c=min\left\{ c_1,c_2\right\}$ and $a_m={3m\choose m}+2m-1$. Let $B_1, B_2, B_3$ and $X$ be disjoint copies of the complete graph, where the cardinalities are defined as follows: $|B_1|=|B_3|={3m\choose m}, |B_2|=3m$ and  $|X|=\frac{1-c}{c} \left(2{3m\choose m}+4m\right)+1$.
The vertex set $X$ is further partitioned into five subsets, denoted by $X_1,X_2,X_3,X_4$ and $X_5$, with the respective cardinalities specified as: $|X_1|=|X_5|=\left(\frac{2-2(c_1+c_2)}{c(2-c_1-c_2)}-1\right)a_m$,
$|X_2|=\frac{2c_1}{c(2-c_1-c_2)}a_m$,
$|X_4|=\frac{2c_2}{c(2-c_1-c_2)}a_m$
and $|X_3|=\frac{2}{c}-1$.

The mixed graph $G_m$ is constructed from $B_1, B_2, B_3$ and $X$ through the following processes:
(1) All directed edges are added from $B_i$ to $X_2$, and from $X_4$ to $B_i$ for $i\in{1,2,3}$;
(2) all undirected edges are added between $B_1\cup B_2$ and $X_1$, as well as between $B_2\cup B_3$ and $X_5$;
(3) for each $2m$-subset $S\subseteq V(B_2)$, one vertex $u_S$ in $B_1$ and one vertex $v_S$ in $B_3$ are associated with $S$,  and undirected edges are then added to joining $u_S$ and $v_S$ to every vertex in $S$.

The degree of each vertex of $G_m$ is presented in Table~$1$. It follows from Table~$1$ that the graph $G_m$ is a $(c_1,c_2)$-mixed graph based on straightforward calculation.

\begin{center}
\small
\begin{tabular}{|c |c |c |c |c |}
\hline &out-degree & in-degree & un-degree & degree \\\hline

$B_1$ & $\frac{2c_1}{c(2-c_1-c_2)}a_m$ & $\frac{2c_2}{c(2-c_1-c_2)}a_m$ & $\frac{2-2(c_1+c_2)}{c(2-c_1-c_2)}a_m$ & $\frac{2}{c(2-c_1-c_2)}a_m$ \\\hline

$B_2$ & $\frac{2c_1}{c(2-c_1-c_2)}a_m$ & $\frac{2c_2}{c(2-c_1-c_2)}a_m$  &$\left(\frac{4-4(c_1+c_2)}{c(2-c_1-c_2)}-1\right)a_m+\frac{1}{3}{3m\choose m}+m$ &$\frac{2-c}{c}a_m+\frac{1}{3}{3m\choose m}+m$  \\\hline

$B_3$ & $\frac{2c_1}{c(2-c_1-c_2)}a_m$ & $\frac{2c_2}{c(2-c_1-c_2)}a_m$ & $\frac{2-2(c_1+c_2)}{c(2-c_1-c_2)}a_m$ & $\frac{2}{c(2-c_1-c_2)}a_m$ \\\hline

$X_1$ & $0$ & $0$ & $\frac{2-c}{c} \left({3m\choose m}+2m\right)+m$ & $\frac{2-c}{c} \left({3m\choose m}+2m\right)+m$  \\\hline

$X_2$ & $0$ & $2{3m\choose m}+3m$ & $\frac{1-c}{c} \left(2{3m\choose m}+4m\right)$ & $\frac{1}{c} \left(2{3m\choose m}+4m\right)-m$ \\\hline

$X_3$ & $0$ & $0$ & $\frac{1-c}{c} \left(2{3m\choose m}+4m\right)$ & $\frac{1-c}{c} \left(2{3m\choose m}+4m\right)$  \\\hline

$X_4$ & $2{3m\choose m}+3m$ & $0$ & $\frac{1-c}{c} \left(2{3m\choose m}+4m\right)$ & $\frac{1}{c} \left(2{3m\choose m}+4m\right)-m$ \\\hline

$X_5$ & $0$ & $0$ & $\frac{2-c}{c} \left({3m\choose m}+2m\right)+m$ & $\frac{2-c}{c} \left({3m\choose m}+2m\right)+m$ \\\hline

\end{tabular}

\vspace{6pt}
Table $1$. $(c_1,c_2)$-mixed graph $G_m$ in Lemma \ref{Lower1}.
\end{center}

Let $D_m$ be an arbitrary orientation of $G_m$. For each $u\in B_1$, we have that $|N_{D_m}(u)\cap V(B_2)|=2m$. Thus, $|N^+_{D_m}(u)\cap V(B_2)| +|N^-_{D_m}(u)\cap V(B_2)|=2m$, where one of these terms cannot exceed $m$. Without loss of generality, we assume $|N^+_{D_m}(u)\cap V(B_2)|\leq m$.
Since $|B_2|=3m$, there exists a subset $R \subseteq V(B_2)-( N^+_{D_m}(u)\cap V(B_2))$ with $|R|=2m$. Consequently, for this subset $R$, there must be a corresponding vertex $v_{R}\in V(B_3)$ such that $D_m$ has no directed $2$-path from $u$ to $v_{R}$ through a vertex of $B_2$. It remains to show that there is also no directed $2$-path from $u$ to $v_{R}$ through a vertex of $X$.

Let $n_m=|V(G_m)|$.
Then $n_m=\frac{1-c}{c}\left(2{3m\choose m}+4m\right)+1+2{3m\choose m}+3m=\frac{2}{c}\left({3m\choose m}+2m\right)-m+1$.
By Table~$1$, we know that the degree of vertices in $B_1\cup B_3$ is the smallest. That is, for each $u\in V(B_1\cup B_3)$,
\begin{align*}
 \delta (G_m)
 &=d_{G_m}(u)\\
 &=\frac{2}{c(2-c_1-c_2)}\left({3m\choose m}+2m-1\right)\\
 &=\frac{1}{2-c_1-c_2}n_m+\frac{1}{2-c_1-c_2}m-
 \frac{c+2}{c(2-c_1-c_2)}.
\end{align*}

Using the Robbins' formula \cite{Robbins}, we get
$${3m\choose m}=\frac{(3m)!}{m!(2m)!}=\frac{(\frac{3m}{e})^{3m}\sqrt{6\pi m}e^{\gamma_1}}{(\frac{m}{e})^{m}\sqrt{2\pi m}e^{\gamma_2}(\frac{2m}{e})^{2m}\sqrt{4\pi m}e^{\gamma_3}}=\frac{\sqrt{3}}{2\sqrt{\pi m}}\left(\frac{27}{4}\right)^m e^{\gamma_1-\gamma_2-\gamma_3},$$
where $\frac{1}{36m+1}<\gamma_1<\frac{1}{36m},
\frac{1}{12m+1}<\gamma_2<\frac{1}{12m},  \frac{1}{24m+1}<\gamma_3<\frac{1}{24m}.$

For sufficiently large $m$, there exist constants $\alpha_1, \beta_1 >0$ satisfying that
$$n_m=\frac{2}{c}\left({3m\choose m}+2m\right)-m+1<\alpha_1{3m\choose m}=\alpha_1 \frac{\sqrt{3}}{2\sqrt{\pi m}}\left(\frac{27}{4}\right)^m e^{\gamma_1-\gamma_2-\gamma_3}
<\beta_1\frac{1}{\sqrt{m}}\left(\frac{27}{4}\right)^{m}.$$

Taking logarithms of both sides, we get
$$m>\frac{\ln n_m+\frac{1}{2}\ln m-\ln\beta_1}{\ln\frac{27}{4}}.$$

Since $m$ is sufficiently large, thus
\begin{align*}
\delta(G_m)&=\frac{1}{2-c_1-c_2}n_m+
\frac{1}{2-c_1-c_2}m-\frac{c+2}{c(2-c_1-c_2)}\\
&>\frac{1}{2-c_1-c_2}n_m+
\frac{1}{2-c_1-c_2} \frac{\ln n_m+\frac{1}{2}\ln m-\ln\beta_1}{\ln\frac{27}{4}}-\frac{c+2}{c(2-c_1-c_2)}\\
&>\frac{1}{2-c_1-c_2}n_m+\frac{1}{2-c_1-c_2}  \frac{\ln n_m}{\ln\frac{27}{4}}.
\end{align*}

Let $n_m=n$, we have
$$\delta(G_m)\geq\frac{1}{2-c_1-c_2}n+\frac{2}{2-c_1-c_2}
  \frac{\ln n}{2\ln \frac{27}{4}}.$$
Thus, $G_m$ is our desired graph.
\end{pf} $\hfill\square$

\begin{lemma}\label{Lower2} Let $c_1$, $c_2$ be two rational numbers with $0\leq c_1, c_2\le 1$. If $\frac{1-c_1-c_2}{c(2-c_1-c_2)}< 1$, then there exists an infinite family of $(c_1,c_2)$-mixed graphs $G_m$ of order $n$ and
$$\delta(G_m)\geq\frac{1}{2-c_1-c_2}n+\frac{2}{2-c_1-c_2}
  \frac{\ln n}{2\ln \frac{27}{4}},$$
whose oriented diameter is at least three.
\end{lemma}
\begin{pf}
Let $c=min\left\{ c_1,c_2\right\}$ and $a_m={3m\choose m}+2m-1$. Let $B_1, B_2, B_3$ and $X$ be disjoint copies of the complete graph, where the cardinalities are defined as follows: $|B_1|= |B_3|={3m\choose m}, |B_2|=3m$ and  $|X|=\frac{1}{1-c_1-c_2}  \left(2{3m\choose m}+4m\right)+1$. The vertex set $X$ is further partitioned into five subsets, denoted by $X_1,X_2,X_3,X_4$ and $X_5$, with the respective cardinalities specified as: $|X_1|=|X_5|=a_m$,
$|X_2|=\frac{2c_1}{1-c_1-c_2}a_m$,
$|X_4|=\frac{2c_2}{1-c_1-c_2}a_m$
and $|X_3|=\frac{2}{1-c_1-c_2}+1$.

We can construct $G_m$ by adding edges similar to the proof of Lemma~\ref{Lower1}.
The degree of each vertex of $G_m$ is shown in Table~$2$. It follows from Table~$2$ that the graph $G_m$ is a $(c_1,c_2)$-mixed graph based on straightforward calculation.

\begin{center}
\small
\begin{tabular}{|c |c |c |c |c |}
\hline &out-degree & in-degree & un-degree & degree \\\hline

$B_1$ & $\frac{2c_1}{1-c_1-c_2}a_m$ & $\frac{2c_2}{1-c_1-c_2}a_m$ & $2a_m$ & $\frac{2}{1-c_1-c_2}a_m$\\\hline

$B_2$ & $\frac{2c_1}{1-c_1-c_2}a_m$ & $\frac{2c_2}{1-c_1-c_2}a_m$ & $3a_m+\frac{1}{3}{3m\choose m}+m$ & $\frac{3-c_1-c_2}{1-c_1-c_2}a_m+\frac{1}{3}{3m\choose m}+m$ \\\hline

$B_3$ & $\frac{2c_1}{1-c_1-c_2}a_m$ & $\frac{2c_2}{1-c_1-c_2}a_m$ & $2a_m$ & $\frac{2}{1-c_1-c_2}a_m$ \\\hline

$X_1$ & $0$ & $0$ & $\frac{3-c_1-c_2}{1-c_1-c_2} \left({3m\choose m}+2m\right)+m$ & $\frac{3-c_1-c_2}{1-c_1-c_2} \left({3m\choose m}+2m\right)+m$ \\\hline

$X_2$ & $0$ & $2{3m\choose m}+3m$ & $\frac{1}{1-c_1-c_2} \left(2{3m\choose m}+4m\right)$ & $\frac{2-c_1-c_2}{1-c_1-c_2} \left(2{3m\choose m}+4m\right)-m$ \\\hline

$X_3$ & $0$ & $0$ & $\frac{1}{1-c_1-c_2} \left(2{3m\choose m}+4m\right)$ & $\frac{1}{1-c_1-c_2} \left(2{3m\choose m}+4m\right)$\\\hline

$X_4$ & $2{3m\choose m}+3m$ & $0$ & $\frac{1}{1-c_1-c_2} \left(2{3m\choose m}+4m\right)$ & $\frac{2-c_1-c_2}{1-c_1-c_2} \left(2{3m\choose m}+4m\right)-m$\\\hline

$X_5$ & $0$ & $0$ & $\frac{3-c_1-c_2}{1-c_1-c_2} \left({3m\choose m}+2m\right)+m$ & $\frac{3-c_1-c_2}{1-c_1-c_2} \left({3m\choose m}+2m\right)+m$ \\\hline

\end{tabular}

\vspace{6pt}
Table $2$. $(c_1,c_2)$-mixed graph $G_m$ in Case~2
\end{center}

Let $D_m$ be an arbitrary orientation of $G_m$. By a similar argument as the proof of Lemma~\ref{Lower1}, $D_m$ has no directed $2$-path from $u$ to $v_{R}$ through a vertex of $B_2$. It remains to show that there is also no directed $2$-path from $u$ to $v_{R}$ through a vertex of $X$.

Let $n_m=|V(G_m)|$.
Then $n_m=\frac{1}{1-c_1-c_2} \left(2{3m\choose m}+4m\right)+1+2{3m\choose m}+3m=\frac{2-c_1-c_2}{1-c_1-c_2}\left(2{3m\choose m}+4m\right)-m+1$.
By simple observation and calculation, we know that the degree of vertices in $B_1\cup B_3$ is smallest. That is, for each $u\in V(B_1\cup B_3)$,
\begin{align*}
 \delta (G_m)
 &=d_{G_m}(u)\\
 &=\frac{2}{1-c_1-c_2}\left({3m\choose m}+2m-1\right)\\
 &=\frac{1}{2-c_1-c_2}n_m+\frac{1}{2-c_1-c_2}m-
 \frac{5-3c_1-3c_2}{(1-c_1-c_2)(2-c_1-c_2)}.
\end{align*}

For sufficiently large $m$, there exist constants $\alpha_2, \beta_2 >0$ satisfying that
$$n_m=\frac{2-c_1-c_2}{1-c_1-c_2}\left(2{3m\choose m}+4m\right)-m+1<\alpha_2{3m\choose m}=\alpha_2\frac{\sqrt{3}}{2\sqrt{\pi m}}\left(\frac{27}{4}\right)^m e^{\gamma_1-\gamma_2-\gamma_3}
<\beta_2\frac{1}{\sqrt{m}}\left(\frac{27}{4}\right)^{m}.$$

Taking logarithms of both sides, we get
$$m>\frac{\ln n_m+\frac{1}{2}\ln m-\ln\beta_2}{\ln\frac{27}{4}}.$$

Since $m$ is sufficiently large, thus
\begin{align*}
\delta(G_m)&=\frac{1}{2-c_1-c_2}n_m+\frac{1}{2-c_1-c_2}m-
 \frac{5-3c_1-3c_2}{(1-c_1-c_2)(2-c_1-c_2)}\\
&>\frac{1}{2-c_1-c_2}n_m+
\frac{1}{2-c_1-c_2} \frac{\ln n_m+\frac{1}{2}\ln m-\ln\beta_2}{\ln\frac{27}{4}}-\frac{5-3c_1-3c_2}{(1-c_1-c_2)(2-c_1-c_2)}\\
&>\frac{1}{2-c_1-c_2}n_m+\frac{1}{2-c_1-c_2}  \frac{\ln n_m}{\ln\frac{27}{4}}.
\end{align*}

Let $n_m=n$, we have
$$\delta(G_m)\geq\frac{1}{2-c_1-c_2}n+\frac{2}{2-c_1-c_2}
  \frac{\ln n}{2\ln \frac{27}{4}}.$$
Thus, $G_m$ is our desired graph.
\end{pf} $\hfill\square$

Follows from Lemmas \ref{Lower1} and \ref{Lower2}, we have the following theorem.

\begin{theorem}\label{Lower} Let $c_1$, $c_2$ be two rational numbers with $0\leq c_1, c_2\le 1$. There exists an infinite family of $(c_1,c_2)$-mixed graphs $G_m$ of order $n$ and
$$\delta(G_m)\geq\frac{1}{2-c_1-c_2}n+\frac{2}{2-c_1-c_2}
  \frac{\ln n}{2\ln \frac{27}{4}},$$
whose oriented diameter is at least three.
\end{theorem}

Combining Theorems \ref{Upper} and \ref{Lower}, we have that $\delta(n,c_1,c_2)=\frac{n}{2-c_1-c_2}+\Theta(\ln n)$.

\begin{remark}
For a undirected graph, that is, $c_1=c_2=0$, we have that $\delta(n,0,0)=\frac{n}{2}+\Theta(\ln n)$. Thus, our results can be viewed as a generalization of the result of the undirected graph by Czabarka, Dankelmann and Sz\'ekely in \cite{Czabarka}.
\end{remark}

\section{Conclusion}

A mixed graph contains both undirected and directed edges, which can be viewed an intermediate stage in the orientation process. To measure the ratio of directed edges in a mixed graph, we introduce the concept of the mixed ratio of a mixed graph. For each mixed graph $G$, there exist rational numbers $c_1$ and $c_2$ with $0\leq c_1,c_2\le 1$ such that $G$ is a {\it $(c_1,c_2)$-mixed graph.

In this paper, we prove} that for a $(c_1,c_2)$-mixed graph $G$ of order $n$, the minimum degree threshold for diameter two orientability is $\frac{n}{2-c_1-c_2}+\Theta(\ln n)$. This work generalizes the results of Czabarka, Dankelmann and Sz\'ekely concerning undirected graphs, since a undirected graph can be view as a $(0,0)$-mixed graph. We conclude our paper by presenting the following two problems.

\begin{problem}
For a $(c_1,c_2)$-mixed graph $G$, determine the exact lower and upper bounds of its oriented radius.
\end{problem}

\begin{problem}
For a $(c_1,c_2)$-mixed graph $G$, determine the exact lower and upper bounds of its oriented diameter.
\end{problem}

For each integer $r\ge 1$, let $g(r,c_1,c_2)$ be the smallest number such that each bridgeless $(c_1,c_2)$-mixed graph with radius $r$ possesses an orientation with a radius of at most $g(r,c_1,c_2)$. For each integer $d\ge 1$, let $f(d,c_1,c_2)$ be the smallest number such that each bridgeless $(c_1,c_2)$-mixed graph with diameter $d$ possesses an orientation with a diameter of at most $f(d,c_1,c_2)$. Following the work of Chv\'atal and Thomassen \cite{Chvatal}, as well as those of Babu, Benson and Rajendraprasad \cite{Babu}, we have the inequalities:
$$r^2+r\le g(r,c_1,c_2)\le 1.5r^2+r+1,$$ and $$\frac{1}{2}d^2+d\le f(d,c_1,c_2)\le 3d^2+2d+2.$$

\section{Acknowledgements}
The paper was supported by Postgraduate Education Reform and Quality Improvement Project of Henan Province (YJS2024KC32).

\end{document}